\documentclass[a4paper,12pt]{article}
\usepackage{a4wide}
\usepackage{amsmath}
\usepackage{amssymb}
\usepackage{amsthm}
\usepackage{latexsym}
\usepackage{graphicx}
\usepackage[english]{babel}
\usepackage{makeidx}
              
\newtheorem{dfn} [section]{Definition}
\newtheorem{obs} [section]{Remark}

\newtheorem{prop}[section]{Proposition}

\newtheorem{teor}[section]{Theorem}
\newtheorem{lema}[section]{Lemma}
\newtheorem{cor} [section]{Corollary}

\begin{document}

\selectlanguage{english}
\frenchspacing

\large
\begin{center}
\textbf{A stable property of Borel type ideals.}

Mircea Cimpoea\c s
\end{center}

\normalsize

\footnotetext[1]{This paper was supported by the CEEX Program of the Romanian
Ministry of Education and Research, Contract CEX05-D11-11/2005 and by 
 the Higher Education Commission of Pakistan.}

\begin{abstract}
In this paper, we extend a result of Eisenbud-Reeves-Totaro in the frame of ideals of Borel type. 

\vspace{5 pt} \noindent \textbf{Keywords:} Borel type ideals, $p$-Borel ideals, Mumford-Castelnuovo regularity.

\vspace{5 pt} \noindent \textbf{2000 Mathematics Subject
Classification:}Primary: 13P10, Secondary: 13E10.
\end{abstract}

\section*{Introduction.}

Let $K$ be an infinite field, and let $S=K[x_1,\ldots,x_n],n\geq 2$ the polynomial ring over $K$.
Bayer and Stillman \cite{BS} note that Borel fixed ideals $I\subset S$ satisfies the following property:
\[(*)\;\;\;\;(I:x_j^\infty)=(I:(x_1,\ldots,x_j)^\infty)\;for\; all\; j=1,\ldots,n.\] Herzog, Popescu and Vladoiu \cite{hpv} define that a monomial ideal $I$ is of \emph{Borel type} if it satisfy $(*)$. We mention that this concept appears also in \cite[Definition 1.3]{CS} as the so called weakly stable ideal. Herzog, Popescu and Vladoiu proved in \cite{hpv} that $I$ is of Borel type, if and only if for any monomial $u\in I$ and for any $1\leq j<i \leq n$ and $q>0$ with $x_i^{q}|u$, there exists an integer $t>0$ such that $x_j^{t}u/x_i^{q}\in I$. As a consequence, it became obvious that the sum of ideals of Borel type is an ideal of Borel type.

We recall that $I_{\geq e}$ is the ideal generated by the monomials of degree $\geq e$ from $I$. We prove that if $I$ is an ideal of Borel type, then $I_{\geq e}$ is stable whenever $e\geq reg(I)$ (Theorem $6$). This allows us to give a generalization of a result of Eisenbud-Reeves-Totaro (Corollary $8$). Also, the result can be easily extended to the monomial ideals $I\subset S$ with $Ass(S/I)$ totally ordered by inclusion (Corollary $9$).

The original motivation of this paper came from the remark that the regularity of the $d$-fixed ideals, see \cite[Definition 1.4]{mir}, and their sums satisfy the upper bound $n \cdot (deg(I)-1)+1$, where $deg(I)$ is the maximum degree of a minimal generator of $I$. 
Later, A.Imran and A.Sarfraz show in \cite{saf} that an ideal of Borel type $I\subset S$ has the regularity upper bounded by $n \cdot (deg(I)-1)+1$. Our idea was to extend a result of Eisenbud-Reeves-Totaro, \cite[Proposition 12]{ert}, in the frame of sums of $d$-fixed ideals, but in fact we realize that it is much simpler to use the more general notion of Borel type ideals.

The author is grateful to his adviser Dorin Popescu for his encouragement and valuable suggestions.
My thanks goes also to the School of Mathematical Sciences, GC University, Lahore, Pakistan for supporting 
and facilitating this research.

\vspace{2mm} \noindent {\footnotesize
\begin{minipage}[b]{15cm}
 Mircea Cimpoea\c s, Institute of Mathematics of the Romanian Academy, Bucharest, Romania\\
 E-mail: mircea.cimpoeas@imar.ro
\end{minipage}}

\newpage
\section*{A stable property of Borel type ideals.}

It would be appropriate to recall the definition of the Castelnuovo-Mumford regularity. We refer the reader to \cite{E} for further details on the
subject.

\begin{dfn}
Let $K$ be an infinite field, and let $S=K[x_1,...,x_n],n\geq 2$ the polynomial ring over $K$.
Let $M$ be a finitely generated graded $S$-module. The \emph{Castelnuovo-Mumford regularity} $reg(M)$ of $M$ is
\[ \max_{i,j} \{j-i :\; \beta_{ij}(M)\neq 0\}.\]
\end{dfn}

\begin{obs}{\em
For any monomial $u\in S$ we denote $m(u)=max\{i:\; x_i|u\}$. For any monomial ideal $I\subset S$, we denote
$m(I)=max\{m(u):\;u\in G(I)\}$, where $G(I)$ is the set of the minimal generators of $I$. Also, if $M$ is a graded $S$-module of finite length, we denote $s(M)=max\{t:\;M_{t}\neq 0\}$.

Let $I\subset S$ be a Borel type ideal. We recall the sequential chain of $I$ as defined in \cite{hpv}:
\[ I = I_{0} \subset I_{1}\subset \cdots \subset I_r = S, \]
where $n_{\ell}=m(I_{\ell})$ and $I_{\ell+1}=(I_{\ell}:x_{n_{\ell}}^{\infty})$. Notice that $r\leq n$, since $n_{\ell}>n_{\ell+1}$ for all $0\leq \ell<r$. 

We define $J_{\ell}\subset S_{\ell} = K[x_1,\ldots,x_{n_{\ell}}]$ to be the ideal generated by $G(I_{\ell})$. We have the following 
formula for the Mumford-Castelnuovo regularity of $I$, \cite[Corollary 2.5]{hpv}, \[(1)\; reg(I) = max\{s(J_{0}^{sat}/J_0),s(J_{1}^{sat}/J_1),\cdots,s(J_{r-1}^{sat}/J_{r-1})\} + 1.\]
Also, \cite[Corollary 2.5]{hpv} states
\[ (2)\; I_{\ell+1}/I_{\ell} \cong (J_{\ell}^{sat}/J_{\ell})[x_{n_{\ell}+1},\ldots,x_n].\]
}
\end{obs}

\begin{obs}
Note that if the length of the sequential chain of $I$ is $r=1$ then $I$ is an Artinian ideal. 
Also, any Artinian monomial ideal is obviously an ideal of Borel type.
\end{obs}

\begin{lema}
Let $I\subset S$ be an ideal and $I'=IS'$ the extension of $I$ in $S'=S[x_{n+1}]$. If $e\geq deg(I)$, then
$I_{\geq e}$ is stable if and only if $I'_{\geq e}$ is stable.
\end{lema}

\begin{proof}
Suppose that $I_{\geq e}$ is stable. Let $u\in I'_{\geq e}$ be a monomial. Then $u=x_{n+1}^{k}\cdot v$ for some monomial $v\in I$. If $k>0$ then $m(u)=n+1$ and therefore, for any $i<n+1$, $x_i\cdot u/x_{n+1} = x_{n+1}^{k-1}\cdot x_i\cdot v \in I'_{\geq e}$. If $k=0$ then
$m(u)\leq n$ and since $I_{\geq e}$ is stable, it follows $x_i \cdot u/ x_{m(u)} \in I'_{\geq e}$. Thus $I'_{\geq e}$
is stable. For the converse, simply notice that $G(I_{\geq e})\subset G(I'_{\geq e})$ and since it is enough to
check the stable property only for the minimal generators, we are done.
\end{proof}

\begin{lema}
If $I\subset S$ is an artinian monomial ideal and $e\geq reg(I)$ then $I_{\geq e}$ is stable.
\end{lema}

\begin{proof}
Note that $reg(I) = s(S/I) + 1$, therefore, if $e\geq reg(I)$ then $I_{\geq e} = S_{\geq e}$, thus $I_{\geq e}$ is stable.
\end{proof}

\begin{teor}
Let $I\subset S$ be a Borel type ideal and let $e\geq reg(I)$ be an integer. Then $I_{\geq e}$ is stable.
\end{teor}

\begin{proof}
We use induction on $r\geq 1$, where $r$ is the length of the sequential chain of $I$. If $r=1$, i.e. $I$ is an
artinian ideal, we are done from the previous lemma. 

Suppose now $r>1$ and let $I = I_0 \subset I_1 \subset \cdots \subset I_r = S$ be the sequential chain of $I$. 
Using the induction hypothesis, we may assume $(I_1)_{\geq e}$ stable for $e\geq reg(I_1)$. On the other hand, from $(1)$ it follows that $reg(I_1)\leq reg(I)$, thus 
$(I_1)_{\geq e}$ is stable for $e\geq reg(I)$.

Since $J_0^{sat} = I_1 \cap S_{n_0}$, using Lemma $4$ it follows that $(J_0^{sat})_{\geq e}$ is stable.
Let $e\geq reg(I)$. Since $reg(I) \geq s(J_0^{sat}/J_0) + 1$ it follows
that $(J_0)_{\geq e} = (J_0^{sat})_{\geq e}$ is stable. Since $I_0 = J_0S$, using again Lemma $4$, we get
$I_{\geq e}$ stable for $e\geq reg(I)$, as required.
\end{proof}

We recall the following result of Eisenbud-Reeves-Totaro.

\begin{prop}\cite[Proposition 12]{ert}
Let $I$ be a monomial ideal with $deg(I)=d$ and let $e\geq d$ such that $I_{\geq e}$ is stable. Then $reg(I)\leq e$.
\end{prop}

\begin{cor}
If $I$ is a Borel type ideal, then $reg(I)=min\{e:\; I_{\geq e}\;is\; stable \}$.
\end{cor}

\begin{proof}
Firstly, note that $reg(I)\geq deg(I)$ since $\beta_{0i}(I)\neq 0$ for any integer $i$ such that there exists a monomial $u\in G(I)$ of degree $i$.
Therefore, the hypothesis of Theorem $6$ and Proposition $7$ are fulfill. Theorem $6$ gives the "$\geq$" inequality and Proposition $7$ gives the other inequality, so we get the required result.
\end{proof}

\begin{cor}
If $I$ is a monomial ideal with $Ass(S/I)$ totally ordered by inclusion, then $reg(I)=min\{e:\; I_{\geq e}\;is\; stable \}$.
\end{cor}

\begin{proof}
By renumbering the variables, we can assume that each $p\in Ass(S/I)$ has the form 
$p=(x_1,x_2,...,x_r)$ for some $1\leq r\leq n$. Therefore, by \cite[Theorem 2.2]{saf}, we may assume
that $I$ is an ideal of Borel type, and then we apply Corollary $8$.
\end{proof}

\begin{cor}
Suppose $I = I_{1} + \cdots + I_{m}$ where $I_{i}$ are Borel type ideals, for $i=1,\ldots,m$. Then
\[ reg(I) \leq max\{reg(I_{1}),\ldots,reg(I_{m})\}.\]
\end{cor}

\begin{proof}
From Corollary $8$ it follows that $(I_{i})_{\geq reg(I_{i})}$ is stable for $i=1,\ldots,m$. So, if $e = max\{reg(I_{1}),\ldots,reg(I_{r})\}$ it follows that $(I_{i})_{\geq e}$ are stables for all $i$. Therefore
$I_{\geq e}=(\sum_{i=1}^{r}I_{i})_{\geq e} = \sum_{i=1}^{r}(I_{i})_{\geq e}$ is stable, as a sum of stable ideals.
Thus, from Proposition $7$ we get $reg(I)\leq e$ as required.
\end{proof}

\end{document}